\documentclass[12pt]{article}
\usepackage{amsmath, amsthm, amsfonts}
\usepackage{amssymb}
\usepackage{amscd}
\usepackage{verbatim}

\begin{document}
\newcommand{\hilb}{{\mathcal{H}_{\mathbf{a}}}}
\newcommand{\merhav}{{\mathcal D}^{1,2}}
\newcommand{\be}{\begin{equation}}
\newcommand{\ee}{\end{equation}}
\newcommand{\bea}{\begin{eqnarray}}
\newcommand{\eea}{\end{eqnarray}}
\newcommand{\bean}{\begin{eqnarray*}}
\newcommand{\eean}{\end{eqnarray*}}
\newcommand{\thkl}{\rule[-.5mm]{.3mm}{3mm}}
\newcommand{\cw}{\stackrel{\rightharpoonup}{\rightharpoonup}}
\newcommand{\id}{\operatorname{id}}
\newcommand{\supp}{\operatorname{supp}}
\newcommand{\wlim}{\mbox{ w-lim }}
\newcommand{\mymu}{{x_N^{-p_*}}}
\newcommand{\R}{{\mathbb R}}
\newcommand{\N}{{\mathbb N}}
\newcommand{\Z}{{\mathbb Z}}
\newcommand{\Q}{{\mathbb Q}}
\newtheorem{theorem}{Theorem}[section]
\newtheorem{corollary}[theorem]{Corollary}
\newtheorem{lemma}[theorem]{Lemma}
\newtheorem{definition}[theorem]{Definition}
\newtheorem{remark}[theorem]{Remark}
\newtheorem{proposition}[theorem]{Proposition}
\newtheorem{conjecture}[theorem]{Conjecture}
\newtheorem{question}[theorem]{Question}
\newtheorem{example}[theorem]{Example}
\newtheorem{Thm}[theorem]{Theorem}
\newtheorem{Lem}[theorem]{Lemma}
\newtheorem{Pro}[theorem]{Proposition}
\newtheorem{Def}[theorem]{Definition}
\newtheorem{Exa}[theorem]{Example}
\newtheorem{Exs}[theorem]{Examples}
\newtheorem{Rems}[theorem]{Remarks}
\newtheorem{Rem}[theorem]{Remark}
\newtheorem{Cor}[theorem]{Corollary}
\newtheorem{Conj}[theorem]{Conjecture}
\newtheorem{Prob}[theorem]{Problem}
\newtheorem{Ques}[theorem]{Question}
\newcommand{\pf}{\noindent \mbox{{\bf Proof}: }}


\renewcommand{\theequation}{\thesection.\arabic{equation}}
\catcode`@=11
\@addtoreset{equation}{section}
\catcode`@=12


\title{A ground state alternative for singular Schr\"odinger operators}
\author{Yehuda Pinchover\\
 {\small Department of Mathematics}\\ {\small  Technion - Israel Institute of Technology}\\
 {\small Haifa 32000, Israel}\\
{\small pincho@techunix.technion.ac.il}\\\and Kyril Tintarev
\\{\small Department of Mathematics}\\{\small Uppsala University}\\
{\small SE-751 06 Uppsala, Sweden}\\{\small
kyril.tintarev@minmail.net}}
\maketitle
\newcommand{\dnorm}[1]{\thkl #1 \thkl\,}

\begin{abstract}
 Let $\mathbf{a}$ be a quadratic form
associated with a Schr\"odinger operator
$L=-\nabla\cdot(A\nabla)+V$ on a domain $\Omega\subset
\mathbb{R}^d$. If $\mathbf{a}$ is nonnegative on
$C_0^{\infty}(\Omega)$, then either there is $W>0$ such that $\int
W|u|^2\,dx\leq \mathbf{a}[u]$ for all
$C_0^{\infty}(\Omega;\mathbb{R})$, or there is a sequence
$\varphi_k\in C_0^{\infty}(\Omega)$ and a function $\varphi>0$
satisfying $L\varphi=0$ such that $\mathbf{a}[\varphi_k]\to 0$,
$\varphi_k\to\varphi$ locally uniformly in
$\Omega\setminus\{x_0\}$. This dichotomy is equivalent to the
dichotomy between $L$ being subcritical resp. critical in
$\Omega$. In the latter case, one has an inequality of Poincar\'e
type: there exists $W>0$ such that for every $\psi\in
C_0^\infty(\Omega;\mathbb{R})$ satisfying $\int \psi \varphi\,dx
\neq 0$ there exists a constant $C>0$ such that $C^{-1}\int
W|u|^2\,dx\le \mathbf{a}[u]+C\left|\int u \psi\,dx\right|^2$ for
all  $u\in C_0^\infty(\Omega;\mathbb{R})$.
\\[2mm]
\noindent  2000 {\em Mathematics Subject Classification.}
Primary 35J10; Secondary  35J20, 35J70, 49R50.\\[1mm] \noindent
{\em Keywords.} Dirichlet form, ground state, quadratic form.
\end{abstract}
\section{Introduction}
Let $\Omega\subset\R^d$ be a domain. We denote $K\Subset \Omega$,
if $K$ is relatively compact in $\Omega$. Let $A:\Omega
\rightarrow \mathbb{R}^{d^2}$ be a measurable matrix valued
function such that for every $K\Subset\Omega$ there exists
$\lambda_K>1$ such that \be \label{stell} \lambda_K^{-1}I_d\le
A(x)\le \lambda_K I_d \qquad \forall x\in K. \ee Let $V\in
L^{p}_{loc}(\Omega;\mathbb{R})$, where $p>\frac{d}{2}$. Throughout
the paper, we  assume that the bilinear form \be \label{a}
\mathbf{a}(u,v):=\int_\Omega\left(A\nabla u\cdot\nabla
\overline{v}+Vu\overline{v}\right)dx \qquad u,v\in
C_0^\infty(\Omega).\ee associated with the Schr\"odinger operator
\be \label{divform} L:=-\nabla\cdot(A\nabla)+V \ee
 is {\em nonnegative} on $C_0^\infty(\Omega)$, that is
\be \label{assume}
\mathbf{a}[u]:=\mathbf{a}(u,u)=\int_\Omega\left(A\nabla u\cdot
\nabla u+V|u|^2\right)dx\ge 0 \qquad \forall u\in
C_0^\infty(\Omega;\mathbb{R}). \ee Under this assumption it is
known \cite{Ag} that the Dirichlet problem is uniquely solvable in
any bounded subdomain. Therefore,  $\mathbf{a}[u]=0$, for $u\in
C_0^\infty(\Omega)$, if and only if $u=0$.  Consequently, the
bilinear form $\mathbf{a}$ defines a scalar product on
$C_0^\infty(\Omega)$. Let $\hilb(\Omega)$ be the completion of
$C_0^\infty(\Omega)$ with respect to the norm
$\sqrt{\mathbf{a}[u]}$.

Generally, $\hilb(\Omega)$ cannot be identified as a space of
measurable functions or even as a space of distributions. More
precisely, it can happen  that there is no continuous imbedding of
$\hilb(\Omega)$ into ${\mathcal D}^\prime(\Omega)$. For example,
if $\Omega=\R^d$, $d=1,2$, and $\mathbf{a}[u]=\int|\nabla
u|^2\,dx$, then $\hilb(\mathbb{R}^d)$ is not a space of
distributions (see e.g. \cite[Section 11.3]{Maz}), and the zero
element of $\hilb(\mathbb{R}^d)$ contains Cauchy sequences in
$C_0^\infty(\mathbb{R}^d)$ that converge in
$L^2_{loc}(\mathbb{R}^d)$ to the constant function. For $d>2$ the
Sobolev inequality
$\left(\int|u|^\frac{2d}{d-2}\,dx\right)^\frac{d-2}{d}\leq
C\int|\nabla u|^2\,dx $ implies that $\hilb(\mathbb{R}^d)$ is
continuously imbedded into $L^\frac{2d}{d-2}(\mathbb{R}^d)$. In
the case $\Omega=\R^d\setminus\{0\}$, $d>2$, and
$\mathbf{a}[u]=\int_{\Omega}(|\nabla
u|^2-(\frac{d-2}{2})^2\frac{|u|^2}{|x|^2})\,dx$, the class of the
zero element of $\hilb(\Omega)$ contains Cauchy sequences in
$C_0^\infty(\Omega)$ that converge in $L^2_{loc}(\Omega)$ to
$C|x|^{-\frac{d-2}{2}}$, with $C\in \mathbb{R}$.

\begin{definition}\label{defnull}{\em A sequence
$\{\varphi_k\}\subset C_0^\infty(\Omega)$ satisfying
$\mathbf{a}[\varphi_k]\to 0$ is called a {\em null sequence} for
the form $\mathbf{a}$. We say that a positive function $\varphi$
is a {\em null state} for the form $\mathbf{a}$, if there exists a
null sequence $\{\varphi_k\}$ for the form $\mathbf{a}$ such that
$\varphi_k\to\varphi$ in $L^2_{loc}(\Omega)$.
 }\end{definition}
\begin{definition}\label{defsgap}{\em
We say that the nonnegative quadratic form $\mathbf{a}$ has a {\em
weighted spectral gap} if there is a function $W>0$ continuous in
$\Omega$ such that
\begin{equation}
\label{truegap} \int_\Omega W|u|^2\,dx\leq \mathbf{a}[u]\qquad
\forall u\in C_0^\infty(\Omega;\mathbb{R}).
\end{equation}
 }\end{definition}
\begin{remark}\label{remnull}{\em
It is easy to see that a null state $\varphi$ is a distributional
solution of the equation $Lu=0$ in $\Omega$. Indeed, suppose that
$L\varphi\neq 0$, then there exists $\psi\in C_0^\infty(\Omega)$
such that $ \int \varphi L \psi\,dx<-\varepsilon<0$. Let
$\varphi_n\in C_0^\infty(\Omega)$ be a null sequence that
converges to $\varphi$ in $L^2_{loc}(\Omega)$. Then
$\mathbf{a}(\varphi_n, \psi)= \int \varphi_n L \psi\,dx<
-\varepsilon/2$ for all $n$ sufficiently large. Consequently,
$\mathbf{a}[\varphi_n+t\psi]\le
\mathbf{a}[\varphi_n]+t^2\mathbf{a}[\psi] -\varepsilon t$, which
is negative for $t>0$ sufficiently small and $n$ sufficiently
large. This contradicts the assumption that $\mathbf{a}\geq 0$. An
alternative proof of this statement can be deduced from the proof
of Lemma \ref{lem43} which in turn shows that a null state is
actually a weak solution.
 }\end{remark}

In the present paper we establish the dichotomy between the
existence of a weighted spectral gap and the existence of a null
state for the form $\mathbf{a}$.

\begin{theorem}
\label{gap} Let $\Omega\subset\R^d$ be a domain, and assume that
the form $\mathbf{a}$ is nonnegative. Then $\mathbf{a}$ has either
a weighted spectral gap or a null state. If $\mathbf{a}$ has a
null state $\varphi$, then there exists a positive $W\in
C(\Omega)$, such that for every $\psi\in
C_0^\infty(\Omega;\mathbb{R})$ satisfying $\int \psi \varphi \,dx
\neq 0$ there exists a constant $C>0$ such that the following
inequality holds:
 \be\label{Poinc}
 C^{-1}\int_{\Omega} W|u|^2\,dx\le a[u]+C\left|\int_{\Omega} \psi
 u\,dx\right|^2\qquad
 \forall u\in C_0^\infty(\Omega;\mathbb{R}).\ee
Moreover, all norms that are induced by the right hand side of
(\ref{Poinc}) with such functions $\psi$ are equivalent.
\end{theorem}

Theorem \ref{gap} answers a natural question that arises in the
context of improving inequalities of the form
\begin{equation}\label{eqhard1}
\int_\Omega V|u^2|\,dx\le C\int_\Omega A\nabla u\cdot\nabla
\bar{u}\,dx,
\end{equation}
where $V\geq 0\,$: What is the difference between a potential $V$
that can be refined, that is, can be replaced in the inequality
(\ref{eqhard1}) by some $V^\prime\gneqq V$, and a potential that
cannot? In other words, what can be said about a nonnegative
Schr\"odinger operator without any weighted spectral gap? The
answer given in Theorem \ref{gap} is of general relevance to the
work on refining the spectral gap inequalities (see for example,
\cite{BFT,BL,BM,D,FT,MarS} and the references therein). As it is
noted in \cite{FT}, attainment in the standard Sobolev space of
the best constant in the Hardy inequality, or in the Hardy
inequality with a further corrected potential cannot serve as a
criterion for the nonexistence of a weighted spectral gap. The
present paper shows that the relevant energy space for the
Schr\"odinger equation is the completion of $C_0^\infty(\Omega)$
with respect to the quadratic form of the operator appended
possibly with a one-dimensional correction. The Poincar\'e type
inequality (\ref{Poinc}) shows that such a completion is
continuously imbedded into some weighted $L^2$-spaces.
 We show that the existence of
a weighted spectral gap is equivalent to the existence of
continuous imbedding of $\hilb(\Omega)$ into the space of
distributions, which in turn implies a continuous imbedding into
$L^p_{loc}(\Omega)$.

Another equivalent formulation of Theorem~\ref{gap} is the
dichotomy between critical and subcritical elliptic operators. Let
$P$ be a linear second-order elliptic operator with real
coefficients which is defined on a domain $\Omega$ (if $P$ is a
symmetric operator of the form (\ref{divform}), then we denote the
operator by $L$). Let $\mathcal{C}_P(\Omega)$ be the cone of all
positive solutions of the equation $Pu=0$ in $\Omega$. Let $W\in
L^p_{loc}(\Omega;\mathbb{R})$, $p>d/2$. The {\em generalized
principal eigenvalue} is defined by
$$\lambda_0(P,W,\Omega):=\sup\{\lambda\in \mathbb{R}\mid
\mathcal{C}_{P-\lambda W}(\Omega)\neq \emptyset \}.$$

Let $K\Subset \Omega$. Recall \cite{Ag,Pinsk} that $u\in
\mathcal{C}_P(\Omega\setminus K)$ is said to be a {\em positive
solution of the operator $P$ of minimal growth in a neighborhood
of infinity in} $\Omega$, if for any $K\Subset K_1 \Subset \Omega$
and any $v\in C(\overline{\Omega\setminus K_1})\cap
\mathcal{C}_P(\Omega\setminus K_1)$, the inequality $u\le v$ on
$\partial K_1$ implies that $u\le v$ in $\Omega\setminus K_1$. A
positive solution $u\in \mathcal{C}_P(\Omega)$ which has minimal
growth in a neighborhood of infinity in $\Omega$ is called a {\em
ground state} of $P$ in $\Omega$.

The operator $P$ is said to be {\em critical} in $\Omega$, if $P$
admits a ground state in $\Omega$. The operator $P$ is called {\em
subcritical} in $\Omega$, if $\mathcal{C}_P(\Omega)\neq
\emptyset$, but $P$ is not critical in $\Omega$. If
$\mathcal{C}_P(\Omega)= \emptyset$, then $P$ is {\em
supercritical} in $\Omega$.

Assume that $P$ is critical in $\Omega'\varsubsetneq \Omega$, then
$P$ is subcritical in any domain $\Omega_1$ such that
$\Omega_1\varsubsetneq \Omega'$, and supercritical in any domain
$\Omega_2$ such that $\Omega'\varsubsetneq \Omega_2\subset
\Omega$. Furthermore, for any nonzero nonnegative function $W$ the
operator $P+W$ is subcritical and $P-W$ is supercritical in
$\Omega'$. Moreover, if $P$ is critical in $\Omega'$, then $\dim
\mathcal{C}_P(\Omega')=1$ (see e.g. \cite{Pinsk}).

If $P$ is subcritical in $\Omega$, then $P$ admits a positive
minimal Green function $G^\Omega_P(x,y)$ in $\Omega$. For each
$y\in \Omega$, the function $G^\Omega_P(\cdot,y)$ is a positive
solution of the equation $Pu=0$ in $\Omega\setminus\{y\}$ that has
minimal growth in a neighborhood of infinity in $\Omega$.
\begin{theorem}
\label{I} Let $\mathbf{a}$ be a nonnegative quadratic form
associated with a Schr\"odinger operator
$L=-\nabla\cdot(A\nabla)+V$ defined on $\Omega$.  Then the
following conditions are equivalent:

(i) The form $\mathbf{a}$ has a weighted spectral gap.

(ii) The space $\hilb(\Omega)$ is continuously imbedded into
$L^2_{loc}(\Omega)$.

(iii) The space $\hilb(\Omega)$ is continuously imbedded into
${\mathcal D}^\prime(\Omega)$.

(iv) The operator $L$ is subcritical in $\Omega$.
\end{theorem}
Since $\mathbf{a}$ is nonnegative, the operator $L$ is either
critical or subcritical in $\Omega$ \cite{Ag}. Thus,
Theorem~\ref{gap} and Theorem~\ref{I} imply:
\begin{corollary}
The form $\mathbf{a}$ has a null state if and only if $L$ is
critical in $\Omega$. In other words, $\varphi$ is a null state
for the form $\mathbf{a}$ if and only if $\varphi$ is a ground
state of the operator $L$.
\end{corollary}

\section{Preliminary results}
Given a domain $\Omega$, we fix an {\em exhaustion}
$\{\Omega_{N}\}_{N=1}^{\infty}$ of $\Omega$, i.e. a sequence of
smooth, relatively compact domains such that $x_0\in \Omega_1\neq
\emptyset$, $\mbox{cl}({\Omega}_{N})\subset \Omega_{N+1}$ and
$\cup_{N=1}^{\infty}\Omega_{N}=\Omega$. Denote by $B(x_0,\delta)$
the ball of radius $\delta$ centered at $x_0$.

Recall that if $P$ is subcritical in $\Omega$, and $W$ has a
compact support in $\Omega$, then there exists $\varepsilon_0>0$
such that $P+\varepsilon W$ is subcritical in $\Omega$ for all
$|\varepsilon|<\varepsilon_0$ (see \cite{P0,Pinsk}). We need the
following stronger assertion.
 \begin{lemma}\label{lemBS}
Suppose that $P$ is a subcritical operator in $\Omega$. Then there
exists a strictly positive function $W$ such that
$\lambda_0(P,W,\Omega)>0$.
\end{lemma}
\pf Let $u\in\mathcal{C}_P(\Omega)$.  It follows from the proof of
Lemma 3.3 in \cite{Pmax}, that it is sufficient  to find a
positive function  $W$ and a constant $C>0$ such that \be
\label{Hbperteq} \int_{\Omega}G^\Omega_{P}(x,z)W(z)u(z)\,dz\leq
Cu(x) \qquad \forall x\in \Omega.\ee Let $\{\chi_N\}$ be a locally
finite partition of unity on $\Omega$ subordinated to the
exhaustion $\{\Omega_N\}$. Since $\chi_N$ has a compact support it
follows from Theorem 2.10 in \cite{P0} and Remark 3.5 in
\cite{Pmax} that there exists $C_N>0$ such that \be
\label{Hbperteq1}
\int_{\Omega}G^\Omega_{P}(x,z)\chi_N(z)u(z)\,dz\leq C_Nu(x) \ee
for all $x\in \Omega$. Let $\varepsilon_N:=2^{-N}/C_N$, and define
$$W(x):=\sum_{N=1}^\infty\varepsilon_n\chi_N(x).$$ Then
\be \label{Hbperteq2}
\int_{\Omega}G^\Omega_{P}(x,z)W(z)u(z)\,dz\leq
 u(x) \ee for all $x\in \Omega$. Hence,
$\lambda_0(P,W,\Omega)>0$.
 \qed
\begin{remark}\label{remBS}{\em
For further necessary and sufficient conditions for $\lambda_0>0$,
see \cite{MS,Pmax} and the references therein.
 }\end{remark}
\begin{corollary}\label{cor1} Suppose that the operator $L$ is of the form
(\ref{divform}), and let $\mathbf{a}$ be the quadratic form
(\ref{a}) associated with $L$. Then

(i) The operator $L$ is subcritical in $\Omega$ if and only if
$\mathbf{a}$ has a weighted spectral gap, i.e. there exists a
strictly positive function $W$ that
 \be\label{e12}\int_\Omega W |u|^2\,dx\leq
 \mathbf{a}[u]\qquad \forall u\in C^\infty_0(\Omega;\mathbb{R}).\ee

(ii) If $L$ is critical in $\Omega$, then for every open set
$B\Subset \Omega$ there is a strictly positive continuous function
$W$ such that
 \be\label{B}\int_\Omega W |u|^2\,dx\leq \mathbf{a}[u]+\int_B
   |u|^2\,dx\qquad \forall u\in C^\infty_0(\Omega;\mathbb{R}).\ee

(iii) The operator $L$ is critical in $\Omega$ if and only if
there is an open set $B_0\Subset\Omega$, such that
$\int_{B_0}|u|^2\,dx$ is not bounded by $\mathbf{a}[u]$ on
$C^\infty_0(\Omega;\mathbb{R})$.

\vskip 2mm

(iv) If the quadratic form $\mathbf{a}$ has a null state, then $L$
is critical in $\Omega$.
\end{corollary}
\begin{proof} Recall that if $P$ is critical in $\Omega$ if and only for any
nonzero nonnegative function $Q$, the operator $P+Q$ is
subcritical in $\Omega$, and $P-Q$ is supercritical in $\Omega$.
It is well known that if $L$ is symmetric, then
$\mathcal{C}_P(\Omega)\neq\emptyset$ if and only if the quadratic
form $\mathbf{a}$ is nonnegative (see \cite{Ag} and the references
therein). Therefore (i)-(iii) follow from Lemma~\ref{lemBS}. Part
(iii) clearly implies (iv).
\end{proof}

The following key lemma is well known in the case of
Schr\"{o}dinger operators (see, e.g. \cite{DS}).
 \begin{lemma}\label{lemqf}  Let $\psi$ be a solution of the equation
$Lu=0$ in a bounded domain $D\in C^1$, and let $v\in
H^1(D;\mathbb{R})$ such that $v\psi=0$ on $\partial D$. Then
\be\label{e1}\int_D A\nabla (v\psi)\cdot\nabla (v\psi) +
V|v\psi|^2dx=\int_D (A\nabla v\cdot\nabla v) |\psi|^2\,dx.\ee
 \end{lemma}
\begin{proof} Apply the Gauss divergence theorem, and calculate.
\end{proof}

\begin{lemma}\label{lem43}
 Suppose that $L$ is critical in $\Omega$, and let $\varphi$ be its ground state.
Let $\{u_N\}\subset C_0^\infty(\Omega)$ be a null sequence, and
assume that $\{u_N\}$ is locally bounded in $L^2$. Then $\{u_N\}$
has a converging subsequence in $L^2_{loc}(\Omega)$ that converges
to  $c\varphi$  with some $c\in \mathbb{R}$. Moreover, any
converging subsequence of $\{u_N\}$ in $L^2_ {loc}(\Omega)$
converges to $c\varphi$ for some $c\in \mathbb{R}$.

If further, for some $B\Subset \Omega$ the sequence $\{u_N\}$ is
normalized such that $\|u_N\|_{L^2(B)}=1$, then $c\neq 0$. Such a
normalized null sequence does exist. In particular, $\varphi$ is a
null state.
\end{lemma}
\begin{proof} For any $K\Subset\Omega$ there exists $C_K>0$
such that \be\label{eq13}(C_K)^{-1}\leq \varphi(x)\leq C_K\qquad
\forall x\in K .\ee Invoking Lemma~\ref{lemqf} with
$v=u_N/\varphi$, and $\psi=\varphi$, and using  (\ref{stell}) and
(\ref{eq13}),  we infer that $\nabla(u_N/\varphi)$ tends to zero
in $L^2_{loc}(\Omega)$. By our assumption  $u_N/\varphi$ is
locally bounded in $L^2(\Omega)$. In light of  the Sobolev compact
embedding in smooth bounded domains, it follows that (up to a
subsequence) $u_N/\varphi$ converges in $L^2_{loc}(\Omega)$.
Therefore, $u_N/\varphi$ converges in $H^1_{loc}(\Omega)$ to a
function $u$ with a zero gradient, consequently,  $u$ is locally
constant in $\Omega$. Since $\Omega$ is connected $u=const.$ in
$\Omega$.
 Hence, (\ref{eq13}) implies that $\{u_N\}$ converges in $L^2_
{loc}(\Omega)$ to $c\varphi$.

From Corollary \ref{cor1} , part (iii), it follows that there
exist $B\Subset\Omega$ and a null sequence $\{u_N\}$ such that
$\|u_N\|_{L^2(B))}=1$. Thus, $\varphi$ is a null state.
\end{proof}

\begin{remark}\label{rem2}{\em
Clearly, without the assumption of locally boundedness in $L^2$ it
is not true that $\{u_N\}$ converges. Take for instance a bounded
smooth domain and $u_N =N\varphi$. }\end{remark}
\begin{remark}\label{rem3}{\em
 In the subcritical case, any null sequence converges to zero in
$L^2_{loc}(\Omega)$.
 }\end{remark}

\begin{lemma}
\label{nondis} If $\mathbf{a}$ admits a null state $\varphi$, then
for any $\psi\in C_0^\infty(\Omega)$ such that $\int\psi
\varphi\,dx\neq 0$,
 the mapping $u\mapsto\int u\psi\,dx$ is not continuous in $\hilb(\Omega)$.
Consequently, $\hilb(\Omega)$ is not continuously imbedded into
${\mathcal D}^\prime(\Omega)$.
\end{lemma}
\begin{proof} Part (iv) of Corollary \ref{cor1} implies that $L$ is critical in $\Omega$.
It follows from Lemma~\ref{lem43} that there exists a null
sequence $\varphi_k\to c\varphi$ in $L^2_{loc}(\Omega)$ with
$\int_B|\varphi_k|^2\,dx=1$, and $c>0$. Thus, $\int
\varphi_k\psi\,dx\to\int c\varphi\psi\,dx\neq 0$, although
$\varphi_k\to 0$ in $\hilb(\Omega)$.
\end{proof}

\section{Poincar\'e inequality and the space $\merhav_{\mathbf{a}}(\Omega)$}

\begin{proof}[{\bf Proof of Theorem~\ref{I}}] The equivalence of (i) and (iv)
follows from Corollary~\ref{cor1}. From (i) follows immediately
(ii) which implies (iii). If (iv) does not hold, then
Lemma~\ref{lem43} and Lemma~\ref{nondis} imply that condition
(iii) is false.
\end{proof}

\begin{proof} [{\bf Proof of Theorem~\ref{gap}}] If the form $\mathbf{a}$ has a weighted
spectral gap, then every null sequence $w_k$ converges to $0$ in
$L^2_{loc}(\Omega)$, so $\mathbf{a}$ has no null state. If the
form $\mathbf{a}$ has no weighted spectral gap, then by
Theorem~\ref{I}, $L$ is critical in $\Omega$, and by
Lemma~\ref{lem43}, the form $\mathbf{a}$ admits a null state.

Let us prove now (\ref{Poinc}). Due to (\ref{B}) it suffices to
verify that for some open $B\Subset\Omega$, \be \int_B
|u|^2\,dx\le C\left(a[u]+ \left|\int_\Omega
u\psi\,dx\right|^2\right). \ee  Assume that this is false. Then
there is a sequence $u_k$ such that $\mathbf{a}[u_k]\to 0$,
$\int_\Omega u_k\psi\,dx\to 0$, and $\int_B |u_k|^2\,dx=1$.  By
(\ref{B}) $u_k$ is bounded in $L^2_{loc}(\Omega)$, so by
Lemma~\ref{lem43}, $u_k\to \lambda \varphi\neq 0$  in
$L^2_{loc}(\Omega)$.  Then $\int_\Omega u_k\psi\,dx\to \lambda
\int_\Omega  \varphi\psi\,dx\neq 0$,  and we arrive at a
contradiction.

Let $\psi_1,\psi_2\in C_0^\infty(\Omega)$ satisfy
 $\int_\Omega \psi_i \varphi\,dx \neq 0$, $i=1,2$. Then the equivalence of norms follows from
\be \left|\int_\Omega u\psi_1\,dx \right|^2\le C
\left(\mathbf{a}[u]+\left|\int_\Omega u\psi_2\,dx
\right|^2\right), \ee which in turn follows from the
Cauchy-Schwartz inequality and (\ref{Poinc}).
\end{proof}

Recall the standard notation $\merhav(\mathbb{R}^d)$ for
$\hilb(\mathbb{R}^d)$, where $d>2$, and
$\mathbf{a}[u]=\int_{\mathbb{R}^d} |\nabla u|^2\,dx$. Therefore,
it is natural to use the notation $\merhav_{\mathbf{a}}(\Omega)$
for the space $\hilb(\Omega)$ in case of a weighted spectral gap,
and for the closure of $C_0^\infty(\Omega)$ in the norm induced by
the right hand side of (\ref{Poinc}) in the case of the existence
of null state.

By analogy with the compactness of local imbeddings in
$\merhav(\R^d)$, we have the following statement.

\begin{proposition}
\label{compact}  The space $\merhav_{\mathbf{a}}(\Omega)$ is
continuously imbedded into $H^1_{loc}(\Omega)$
 (and therefore, it is compactly imbedded into
 $L^2_{loc}(\Omega)$).
\end{proposition}
\begin{proof}
Consider a ball $\mathcal{B}\subset\merhav_{\mathbf{a}}(\Omega)$.
Fix $\psi\in \mathcal{C}_P(\Omega)$. From (\ref{e1}) it follows
that \be \label{Vpsi} \int_\Omega(A\nabla  (u/\psi)\cdot\nabla
(u/\psi)) |\psi|^2\,dx=\mathbf{a}[u] \qquad \forall u\in
C_0^\infty(\Omega;\mathbb{R}).\ee By density, this implies that
the set $\mathcal{B}_\psi=\{u/\psi:\;u\in \mathcal{B}\}$ is
locally bounded with respect to the Dirichlet norm (i.e. it is
bounded in $\merhav_{loc}(\Omega)$). At the same time, from either
(\ref{truegap}) or (\ref{Poinc}), it follows that $\mathcal{B}$ is
bounded in $L^2_{loc}(\Omega)$. Consequently, using the Leibniz
product rule and the Young inequality, we infer that $\mathcal{B}$
is bounded in $\merhav_{loc}(\Omega)$, and thus also in
$H^1_{loc}(\Omega)$.
\end{proof}

\section{Null sequence converging locally uniformly to null state}
Let $P$ be a second-order elliptic operator with real coefficients
which is not necessarily symmetric and which is defined on a
domain $\Omega\subset \R^d$. Assume that $P$ is critical in
$\Omega$. Fix $x_0\in \Omega$, and let $\varphi$ be the ground
state, satisfying $\varphi(x_0)=1$. Let
$\{\Omega_{N}\}_{N=1}^{\infty}$ be an {\em exhaustion} of
$\Omega$. Without loss of generality assume that $x_0=0$ and
$B(0,1)\Subset \Omega_1$.

We begin this section with the following lemma
 (cf. \cite[Theorem 1.2]{P2}).

 \begin{lemma}\label{lem1}
 Suppose
that $P$ is critical in $\Omega\subset \R^d$. Let $x_1\in
\Omega_1$,  $0<|x_1|<1$.  Consider the function
 $$\psi_N(x):=
  \begin{cases}
    \frac{G_P^{\Omega_N}(x,0)}{G_P^{\Omega_N}(x_1,0)} & x\in\Omega_N,\; x\neq 0,  \\
    \frac{1}{\varphi(x_1)} & x=0.
  \end{cases}
 $$
Then $$\lim_{N\to
\infty}\psi_N(x)=\frac{\varphi(x)}{\varphi(x_1)}\,,$$ locally
uniformly in $\Omega\setminus \{0\}$.
\end{lemma}
\pf  By criticality,
$$\lim_{N\to
\infty} G_P^{\Omega_N}(x_1,0)= \infty.$$  Therefore, \be\label{e2}
\lim_{N\to
\infty}\left(\frac{G_P^{\Omega_1}(x,0)}{G_P^{\Omega_N}(x_1,0)}\right)=0,\ee
locally uniformly in $B(0,|x_1|)\setminus \{0\}$.

Consider the function
 \be\label{e5}\widetilde{\psi}_N(x):=
 \frac{G_P^{\Omega_N}(x,0)-G_P^{\Omega_1}(x,0)}{G_P^{\Omega_N}(x_1,0)}\,.\ee
Note that the function $\widetilde{\psi}_N$ has a removable
singularity at the origin, therefore it may be considered as  a
positive solution of the equation $Pu=0$ in $\Omega_1$. Moreover
by (\ref{e2}), $\widetilde{\psi}_N(x_1)= 1+o(1)$.
 Therefore, by Harnack's inequality
  $c^{-1}\leq\widetilde{\psi}_N(0)\leq c$.  By a standard elliptic argument
$\{\widetilde{\psi}_N\}$ has a subsequence
$\{\widetilde{\psi}_{N_k}\}$ which converges uniformly in any
compact set $K\Subset\Omega_1$, to a positive solution $u$
satisfying $u(x_1)=1$. Hence $\{\psi_{N_k}\}$ converges locally
uniformly in any punctured ball $B(0,r)\setminus
\{0\}\Subset\Omega_1$ to $u$.

On the other hand, $\psi_{N_k}(x)$ is a  positive solution of the
equation $Pu=0$ in $\Omega_{N_k}\setminus \{0\}$, and
$\psi_{N_k}(x_1)=1$. Therefore, the sequence $\{\psi_{N_k}(x)\}$
has a subsequence which converges locally uniformly in
$\Omega\setminus \{0\}$ to a positive solution $u_1(x)$ of the
equation $Pu=0$ in $\Omega\setminus \{0\}$ and $u_1(x_1)= 1$. It
follows that in $\Omega_1\setminus \{0\}$ we have $u_1=u$, and
therefore this subsequence converges to a global positive solution
$u_1$ of the equation $Pu=0$ in $\Omega$, and by uniqueness,
$u_1(x)=\frac{\varphi(x)}{\varphi(x_1)}$. Moreover, since this is
true for any subsequence, it follows that
$$\lim_{N\to \infty}\psi_N(x)=\frac{\varphi(x)}{\varphi(x_1)}$$
locally uniformly in $\Omega\setminus \{0\}$.
 \qed

We use the following cutoff functions.  For $d\geq 3$ define \bean
a_N(x)=
  \begin{cases}
  1 & |x|>\frac{2}{N}, \\
  N(|x|-\frac{1}{N}) & \frac{1}{N}\leq |x|\leq \frac{2}{N},\\
  0& |x|<\frac{1}{N}\,.
  \end{cases}
 \eean
For $d=2$, and for $M<N$, define \bean a_{N,M}(x)=
  \begin{cases}
  1 & |x|>\frac{1}{M}, \\
  \frac{\log |x|N}{\log\frac{N}{M}} & \frac{1}{N}\leq |x|\leq \frac{1}{M},\\
  0& |x|<\frac{1}{N}\,,
  \end{cases}
 \eean
 and denote $a_N:=a_{N,\sqrt{N}}\,$. We have
 \begin{theorem}\label{thm1}
Suppose that the operator $L$ of the form (\ref{divform}) is
critical in $\Omega$, and let $\varphi$ be its ground state. There
exists a null sequence $\{u_N\}\subset H^1(\Omega)$  such that
$\mbox{{\rm supp\,}}u_N\subset\Omega_N$ and $\{u_N\}$ converges
locally uniformly in $\Omega\setminus \{0\}$ to $\varphi$.
\end{theorem}
\pf Set $C:=\max_{|x|\leq 1} \varphi(x)$. By Lemma~\ref{lem1}, the
sequence $\{\varphi(x_1)\psi_{N}\}$ converges locally uniformly in
$\Omega\setminus \{0\}$ to $\varphi$. Hence, there exists an
increasing subsequence $\{M_N\}_{N=1}^\infty \subset \mathbb{N}$
such that
$$\sup_{\frac{1}{N}\leq |x|\leq 1} \varphi(x_1)\psi_{M_N}(x)\leq 2C.$$
Consider the function  $u_N(x):=a_N(x) \varphi(x_1)\psi_{M_N}(x)$.
It follows that  $\{u_N\}$ converges locally uniformly in
$\Omega\setminus \{0\}$ to $\varphi$.

Note that \be\label{intcutoff}\lim_{N\to\infty}\int_\Omega|\nabla
a_N|^2\,dx=0.\ee On the other hand, by the definition of $M_N$,
$0<\varphi(x_1)\psi_{M_N}(x)\leq 2C$ for all $\frac{1}{N}\leq
|x|\leq 1$ and $N\geq 1$. Therefore, (\ref{intcutoff}) and
(\ref{stell}) imply that
\be\label{intcutoff1}\lim_{N\to\infty}\int_\Omega A\nabla
a_N\cdot\nabla a_N (\varphi(x_1)\psi_{M_N}(x))^2\,dx=0.\ee
 Now, use (\ref{e1}) with $v=a_N$, and
$\psi=\varphi(x_1)\psi_{M_N}(x)$, and (\ref{intcutoff1}) to verify
that $\{u_N\}$ is indeed a null sequence.
 \qed

\section{Finding a critical potential for a subcritical operator}
If $P$ is subcritical, then the set $S_+$ (Resp., $S$) of all
continuous $W$, such that $P+W$ is subcritical (Resp., $P+W$ is
not supercritical) is convex \cite{P2}, and as we have seen,
contains a positive function. The set $S_0$ of all continuous $W$,
such that $P+W$ is critical is contained in the set of all extreme
points of $S$. Furthermore, it is known that if $P$ is subcritical
and $W$ is a continuous function that takes a positive value at
some point, and decays sufficiently fast at infinity (for example
$W$ has a compact support), then there is $t_0$ such that $P-tW$
is subcritical for all $0 \leq t<t_0$, critical for $t=t_0$ and
supercritical for $t>t_0$. The precise sufficient decay condition
for $W$ to have such a property depends on $P,W$ and $\Omega$ via
the Green function. An almost optimal general condition is that W
is a small or even semismall perturbation (notions that were
introduced by Y.~Pinchover and M.~Murata, respectively, see, for
example \cite{Murata}).

A critical positive weight $W$ can be found, in particular, via a
minimizer of a nonlinear problem.
\begin{theorem}
\label{shrp} Let $\mathbf{a}$ be the quadratic form associated
with a symmetric subcritical operator $L$ in $\Omega$. Let $W_0\in
L^\infty_{loc}(\Omega)$ be a nonnegative nonzero function, and let
$p> 2$. If the minimum of the following constrained problem
 \be
\label{varst} \kappa=\inf\left\{\mathbf{a}[u] \mid
u\in\hilb(\Omega),\;\int_\Omega W_0|u|^p\,dx=1\right\} \ee is
attained at some positive $v\in\hilb(\Omega)$, then the operator
$L-W$ with $W=\kappa W_0v^{p-2}$ is critical in $\Omega$.
\end{theorem}
Note that the conditions of the theorem fail if $\hilb$ is not
continuously imbedded into $L^p$, in particular, if $L$ is
critical or $p>\frac{d+2}{d-2}$. Moreover, these conditions
imply that $\kappa>0$.
\begin{proof}
If the minimum in (\ref{varst}) is attained at $v>0$, then
$Lv=\kappa W_0v^{p-1}$. In other words,
$v\in\mathcal{C}_{L-W}(\Omega)$. In particular,  the operator
$L-W$ is not supercritical. Let us show that $v$ is a null state
for the quadratic form $\mathbf{b}$ associated with $L-W$. Since
$v\in\hilb(\Omega)$, there exists $v_k\in C_0^\infty$ such that
$v_k\to v$ in $\hilb(\Omega)$. Theorem \ref{I} implies that
$v_k\to v$ in $L^2_{loc}(\Omega)$. By the Fatou lemma
$\liminf\int_\Omega W_0v^{p-2}|v_k|^2\,dx\ge 1$. Thus,
$0\le\limsup\mathbf{b}[v_k]=\lim\mathbf{a}[v_k]-\liminf\int_\Omega
W|v_k|^2\,dx\le \kappa-\kappa=0$. Hence, $v$ is a null state.
\end{proof}

The assumption that the minimum of problem (\ref{varst}) is
attained is not trivial, and its verification typically requires a
concentration-compactness argument. In a recent work \cite{TT}, it
is shown that the conditions of Theorem~\ref{shrp} are satisfied
when $\Omega=\R^d\setminus\R^m$, $1\le m\le d-2$,
$p=\frac{2d}{d-2}\,$, $W_0=1$, and $L=-\Delta -
(\frac{m-2}{2})^2\rho(x)^{-2}$, where $\rho$ is the distance
function to $\partial \Omega$. Consequently, the minimum point $v$
for \be \kappa
=\inf_{\int_{\Omega}|u|^\frac{2d}{d-2}\,dx=1}\;\int_\Omega
\left(|\nabla
u|^2-\left(\frac{m-2}{2}\right)^2\frac{|u|^2}{\rho(x)^{2}}\right)\,dx
\ee is attained, and the operator
$-\Delta-(\frac{m-2}{2})\rho(x)^{-2}-\kappa v^\frac{4}{d-2}$ is
critical. \vskip 1cm
\begin{center}
{\bf Acknowledgments} \end{center} The authors wish to thank
S.~Agmon, S.~Filippas, V.~Maz'ya, and A.~Tertikas for valuable
discussions. This research was partly done at the Technion as
K.~T. was a Lady Davis Visiting Professor, and continued at the
University of Queensland, where K.~T. was supported by the Ethel
Raybould Visiting Fellowship, and a grant from the Swedish
Research Council. The work of Y.~P. was partially supported by the
Fund for the Promotion of Research at the Technion.

\end{document}